%%%%%%%%%%%%%%%%%%%  BEGINNING OF pure.tex  %%%%%%%%%%%
%%%%%%%%%% updated 24 May 1997 
%
\input amstex
\documentstyle{amsppt}
\loadbold
\def\cstar{$C^*$-algebra}
\def\esg{$E_0$-semigroup}

\def\<{\left<}										%for inner products
\def\>{\right>}
\define\tr{\text{trace}}

\magnification=\magstep 1

\topmatter
\title Pure $E_0$-semigroups and absorbing states
\endtitle

\author William Arveson
\endauthor

\affil Department of Mathematics\\
University of California\\Berkeley CA 94720, USA
\endaffil

\date 15 September 1996
\enddate
\thanks This research was supported by
NSF grant DMS95-00291
\endthanks
\keywords von Neumann algebras, automorphism groups, 
\esg s, minimal dilations, completely positive maps
\endkeywords
\subjclass
Primary 46L40; Secondary 81E05
\endsubjclass
\abstract 
An \esg\ $\alpha = \{\alpha_t: t\geq 0\}$ acting on 
$\Cal B(H)$ is called {\it pure} if its tail von Neumann
algebra is trivial in the sense that 
$$
\cap_t\alpha_t(\Cal B(H)) = \Bbb C\bold 1.  
$$
We determine all pure \esg s which have a {\it weakly continuous}
invariant state $\omega$ and which are minimal in an appropriate 
sense.  In such cases the dynamics of the state space must 
stabilize as follows: for every normal state 
$\rho$ of $\Cal B(H)$ there is convergence to equilibrium
in the trace norm 
$$
\lim_{t\to\infty}\|\rho\circ\alpha_t-\omega\|=0.  
$$
A normal state $\omega$ with this property is 
called an {\it absorbing} state for $\alpha$.  

Such \esg s must be cocycle perturbations of 
$CAR/CCR$ flows, and we develop systematic 
methods for constructing those perturbations
which have absorbing states with prescribed 
finite eigenvalue lists.  
\endabstract

\endtopmatter
\vfill\eject
%Replace \pagebreak below with the line above
%to fill the lower part of the title page with
%space, rather than stretching it.
%\pagebreak

\document

\subheading{Introduction}

An \esg\ is a semigroup of normal $*$-endomorphisms 
$\alpha=\{\alpha_t: t\geq 0\}$ of the algebra $\Cal B(H)$ 
of all bounded operators on a separable Hilbert space, which 
satisfies $\alpha_t(\bold 1)=\bold 1$ and the 
natural continuity property
$$
\lim_{t\to 0}\<\alpha_t(x)\xi,\eta\>=\<x\xi,\eta\>, 
\qquad x\in\Cal B(H),\quad\xi,\eta\in H.  
$$ 
There is a sequence of \esg s 
$\alpha^n$, $n=1,2,\dots,\infty$ that can be constructed 
using the {\it natural} irreducible representations of 
either the canonical anticommutation relations or 
the canonical commutation relations. These 
\esg s are called $CAR/CCR$ flows.  They 
occupy a position in 
the category of \esg s roughly analogous to that of the 
unilateral shifts (of various multiplicities) in the 
category of isometries on Hilbert space.  

This paper addresses the perturbation 
theory of $CAR/CCR$ flows. 
A {\it cocycle perturbation} of an \esg\ $\alpha$ is 
an \esg\ $\beta$ which is related to $\alpha$ by way of
$$
\beta_t(x)=U_t\alpha_t(x)U_t^*, \qquad x\in\Cal B(H), t\geq 0
$$
where $\{U_t: t\geq 0\}$ is a strongly continuous family 
of unitary operators in $\Cal B(H)$ which satisfies 
the cocycle equation 
$$
U_{s+t}=U_s\alpha_s(U_t), \qquad s,t\geq 0.  
$$
We are interested in  cocycle perturbations $\beta$ of 
the $CAR/CCR$ flows whose dynamics ``stabilize" in 
that there should exist a normal state $\omega$ 
which is absorbing in the sense that for every 
normal state $\rho$ we have
$$
\lim_{t\to\infty}\|\rho\circ\beta_t - \omega\|=0.   
 \tag{0.1}
$$
It is obvious that 
when an absorbing state exists it is invariant under 
the action of $\beta$, and is in fact the unique normal 
$\beta$-invariant state.  Physicists refer to the property
(0.1) as {\it return to equilibrium}, while in 
ergodic theory the corresponding property is 
called {\it mixing}.  

Every normal state $\omega$ of $\Cal B(H)$ has a unique  
{\it eigenvalue list}, that is, a finite or infinite 
sequence of positive numbers $\lambda_1,\lambda_2, \dots$
which is decreasing ($\lambda_k\geq\lambda_{k+1}$, 
$k\geq 1$) and which has the property 
that for some orthonormal set $\xi_1,\xi_2,\dots$ 
in $H$ we have 
$$
\omega(x)=\sum_k\lambda_k\<x\xi_k,\xi_k\>.  
$$
Clearly $\lambda_1+\lambda_2+\dots = 1$, and of course there 
may be a finite number of repetitions of a given element
in the eigenvalue list.  The {\it set} 
$\{\lambda_k: k\geq 1\}\cup\{0\}$
determined by the eigenvalue list is the 
spectrum of the density operator of $\omega$.  
The eigenvalue list is finite iff $\omega$ is continuous in 
the weak operator topology of $\Cal B(H)$.  

If $\beta$ has an absorbing state $\omega$ then it is obvious from 
(0.1) that the eigenvalue list of $\omega$ contains all of the 
information that could be obtained from the 
dynamics of expectation values 
observed over the long term.  Thus it is natural to ask what 
the possibilities are, and how one finds absorbing states 
for cocycle perturbations of the simplest \esg s.  

In this paper we will be concerned with 
{\it pure} \esg s, i.e., \esg s $\beta$ with the property
that the tail von Neumann algebra is trivial,
$$
\cap_t \beta_t(\Cal B(H))=\Bbb C\bold 1.  \tag{0.2}
$$  
After discussing 
the relationship between purity and the existence of 
absorbing states in general, we take up the analysis of 
{\it weakly continuous} absorbing states, and  
we obtain more or less complete information 
about how to construct them.  Those results are
applied in section 5 to establish the following

\proclaim{Theorem A}
Let $\alpha^n$ be the $CAR/CCR$ flow of 
index $n$, $1\leq n\leq \infty$, and let 
$\lambda_1,\dots,\lambda_r$ be a finite decreasing 
sequence of positive numbers summing to $1$.  
Then there is a cocycle perturbation $\beta$ 
of $\alpha^n$ which has an absorbing state 
$\omega$ with eigenvalue list 
$\lambda_1,\dots,\lambda_r$.  

If  $n\leq r^2-1$ (and in this event $r\geq 2$)
then one can arrange that $\beta$ 
is minimal over the support projection of $\omega$.  

Conversely, if $r\geq 2$ and  $\beta$ is any \esg\ 
which has an 
absorbing state $\omega$ with eigenvalue list 
$\lambda_1,\dots,\lambda_r$, and which is minimal 
over the support projection of $\omega$, then 
$\beta$ is conjugate to 
a cocycle perturbation of $\alpha^n$ for some
$n$, $1\leq n\leq r^2-1$.  
\endproclaim
\remark{Remarks}
The assertions about minimality relate to dilation theory.  
If $\omega$ is an invariant normal state for an \esg\ 
$\beta$ then the support projection $p$ of $\omega$ is 
{\it increasing} in the sense that 
$$
\beta_t(p)\geq p, \qquad t\geq 0,  
$$
(see the discussion following Proposition 2.4).  
It follows that the family of 
completely positive linear maps
$P=\{P_t: t\geq 0\}$ defined on the 
hereditary subalgebra $p\Cal B(H)p\cong \Cal B(pH)$ by
$$
P_t(x)=p\beta_t(x)p,\qquad x\in p\Cal B(H)p, t\geq 0
$$
is in fact a semigroup of completely positive maps.  The 
minimality assertions of the second and third paragraphs
mean that 
$\beta$ is a minimal dilation of $P$ 
in the sense of \cite{2}.  

If $\beta$ is not a minimal dilation of $P$ then there 
is a projection $q\geq p$ satisfying $\beta_t(q)=q$ for 
every $t\geq 0$ and such that the compression of 
$\beta$ to the hereditary subalgebra defined by $q$ 
is a minimal dilation of $P$ (see \cite{2}).  
Thus we may conclude that
\esg s having absorbing states with {\it finite} eigenvalue
lists $\lambda_1,\dots,\lambda_r$, $r\geq 2$ 
are always associated with perturbations of 
$CAR/CCR$ flows.  
\endremark

\remark{Remarks}
In \cite{12}, Powers constructed a new class of examples
of \esg s.  Such an \esg \ $\alpha$ has the property (0.2)
and moreover, there is a unit vector $\xi\in H$ such that the 
pure state $\omega(x) = \<x\xi, \xi\>$ is invariant under 
the action of $\alpha$; indeed $\omega$ is an absorbing 
state.  

In \cite{9}, Bratteli, Jorgensen and Price took up the 
construction of pure invariant states for single 
endomorphisms $\alpha$ of $\Cal B(H)$ satisfying the
discrete counterpart of (0.2)
$$
\cap_n\alpha^n(\Cal B(H)) = \Bbb C\cdot\bold 1,
$$
and they obtain a (non-smooth) paramaterization of 
such states.  While both of these results clearly bear 
some relation to the problems taken up below, we are 
concerned here with absorbing states that are {\it not}
pure.  Indeed, Theorem A has little content for 
eigenvalue lists of length $1$, and the dilation theory 
associated with a pure invariant state is trivial.  

Finally, it is appropriate to comment briefly on 
terminology.  A semigroup of isometries 
$U = \{U_t: t\geq 0\}$ acting on a Hilbert space
$H$ is traditionally called {\it pure} if 
$$
\cap_{t>0}U_tH = \{0\}.  
$$
A familiar theorem in operator theory asserts that every 
pure semigroup of isometries is unitarily equivalent 
to a direct sum of copies of the {\it shift} 
semigroup $S=\{S_t : t\geq 0\}$, which acts on 
the Hilbert space $L^2[0,\infty)$ by way of 
$$
S_tf(x) = 
\cases 
f(x-t),\qquad &x>t\\
0,\qquad &0\leq x\leq t. 
\endcases 
$$
In the theory of \esg s, the proper analogue of 
the {\it shift} of mulitipicity $n=1,2,\dots,\infty$ 
is the $CAR/CCR$ flow of index $n$.  There is no 
theorem in \esg\ theory analogous to the one 
cited above for semigroups of isometries.  Indeed, 
the work of Powers \cite{12}, \cite{13}
implies that
there are \esg s $\alpha$ having the property (0.2)
which are not cocycle conjugate to 
$CAR/CCR$ flows.  Thus we have elected to use the
term {\it pure} for an \esg\ satisfying the 
condition (0.2), and we reserve the term {\it shift}
for the $CAR/CCR$ flows.   
\endremark

\subheading{1.  Purity and absorbing states}

In this section we collect some basic observations about 
pure \esg s acting on von Neumann algebras.  An \esg\ 
$\alpha=\{\alpha_t: t\geq 0\}$ 
acting on a von Neumann algebra $M$ is called {\it pure}
if the intersection $\cap_t\alpha_t(M)$ reduces to the 
scalar multiples of the identity.  The following
result characterizes purity in terms of the action of $\alpha$
on the predual of $M$.  

\proclaim{Proposition 1.1}
Let $\alpha = \{\alpha_t: t\geq 0\}$ be an 
\esg\ acting on a von Neumann algebra $M$.  Then 
$\cap_t\alpha_t(M)=\Bbb C\bold 1$ 
iff for every pair of normal states $\rho_1, \rho_2$ 
of $M$ we have 
$$
\lim_{t\to\infty}\|\rho_1\circ\alpha_t - \rho_2\circ\alpha_t\|=0.  
\tag{1.1.1}
$$
\endproclaim
\demo{proof} We write $M_\infty$ for the von Neumann 
subalgebra $\cap_t\alpha_t(M)$.  
Assume first that $\alpha$ satisfies condition (1.1.1).  To 
show that $M_\infty\subseteq\Bbb C\bold 1$ 
it suffices to show that
for every normal linear functional $\lambda\in M_*$ satisfying 
$\lambda(\bold 1) =0$, we have $\lambda(M_\infty)=\{0\}$.  
Choose such a $\lambda$ and let $\lambda = \lambda_1+i\lambda_2$ 
be its Cartesian decomposition, where 
$\lambda_k(z^*) = \bar\lambda_k(z)$, $k=1,2$.  Since 
$\lambda_k(\bold 1)=0$, it suffices to prove the assertion 
for self-adjoint elements $\lambda$ in the predual of $M$.  

Now by the Hahn decomposition, every self-adjoint element of the 
predual of $M$ which annihilates
the identity operator is a scalar multiple of 
the difference of two normal states.  
Thus, after rescaling, we can assume that there are normal states 
$\rho_1,\rho_2$ of $M$ such that 
$\lambda = \rho_1-\rho_2$, and have to show that 
$\rho_1(x)=\rho_2(x)$ for every element $x\in M_\infty$.  
Since the restriction of each $\alpha_t$ to $M_\infty$ is 
obviously a $*$-automorphism of $M_\infty$, we can find 
a family of operators $x_t\in M_\infty$ such that 
$\alpha_t(x_t) = x$ for every $t\geq 0$.  We have 
$\|x_t\| = \|\alpha_t(x_t)\| = \|x\|$ for every $t$ and
hence 
$$
|\rho_1(x)-\rho_2(x)| = 
|(\rho_1\circ\alpha_t-\rho_2\circ\alpha_t)(x_t)| \leq
\|\rho_1\circ\alpha_t-\rho_2\circ\alpha_t\|\cdot \|x\|
$$
for every $t$.  By hypothesis the right side tends to $0$ with $t$, 
and we have the desired conclusion $\lambda(x)=\rho_1(x)-\rho_2(x)=0$.    

For the converse, let $\rho$ be an arbitrary normal linear 
functional on $M$.  We claim that 
$$
\lim_{t\to\infty}\|\rho\circ\alpha_t\| = 
\|\rho\restriction_{M_\infty}\|.  \tag{1.2}
$$ 
For this, we note first that 
$$
\|\rho\circ\alpha_t\|=\|\rho\restriction_{\alpha_t(M)}\|.  \tag{1.3}
$$
Indeed, the inequality $\leq$ follows from the fact that 
for every $x\in M$, 
$$
|\rho(\alpha_t(x))|\leq 
\|\rho\restriction_{\alpha_t(M)}\|\cdot\|\alpha_t(x)\|.
$$
While on the other hand, if $x\in\alpha_t(M)$ is an 
element of norm $1$ for which 
$$
|\rho(x)|=\|\rho\restriction_{\alpha_t(M)}\|
$$
then we may find $x_0\in M$ with $x=\alpha_t(x_0)$. Noting that
$\|x_0\|=\|x\|$ because $\alpha_t$ is an isometry, we have 
$$
\|\rho\restriction_{\alpha_t(M)}\|=|\rho(x)| =
|\rho\circ\alpha_t(x_0)|\leq \|\rho\circ\alpha_t\|.  
$$

Thus, (1.2) is equivalent to the assertion
$$
\lim_{t\to\infty}\|\rho\restriction_{\alpha_t(M)}\| =
\|\rho\restriction_{M_\infty}\|.  \tag{1.4}
$$
Since the range of $\alpha_t$ is a von Neumann subalgebra
of $M$, we may deduce (1.4) from general principles.  Indeed,
if $M_t$, $t\geq 0$ is a decreasing family of weak$^*$-closed
linear subspaces of the dual of a Banach space $E$ 
having intersection $M_\infty$,  
and $\rho$ is a weak$^*$-continuous linear functional on 
$E^\prime$, then by a standard argument using 
weak$^*$-compactness of the unit ball of $E^\prime$ we find 
that the norms $\|\rho\restriction_{M_t}\|$ must decrease to 
$\|\rho\restriction_{M_\infty}\|$.  

Assuming now that $M_\infty=\Bbb C\bold 1$, let 
$\rho_1$ and $\rho_2$ be normal states of $M$ and let
$\lambda=\rho_1-\rho_2$.  Then the restriction of 
$\lambda$ to $M_\infty$ vanishes, so by (1.2) we have
$$
\lim_{t\to\infty}\|\rho_1\circ\alpha_t-\rho_2\circ\alpha_t\| =
\lim_{t\to\infty}\|\lambda\circ\alpha_t\| = 0, 
$$
as required. \qed
\enddemo

\proclaim{Definition 1.5}
Let $\alpha=\{\alpha_t: t\geq 0\}$ be an \esg\ acting on
a von Neumann algebra $M$.  An  absorbing state for 
$\alpha$ is a normal state $\omega$ on $M$ such that for every
normal state $\rho$, 
$$
\lim_{t\to\infty}\|\rho\circ\alpha_t-\omega\|=0.  
$$
\endproclaim

\remark{Remarks}
An absorbing state $\omega$ is obviously {\it invariant} in the sense 
that $\omega\circ\alpha_t=\omega$, $t\geq 0$, 
and in fact is the {\it unique} 
normal invariant state.  Pure absorbing states for \esg s acting on 
$\Cal B(H)$ were introduced by Powers \cite{13} in his work in 
\esg s of type $II$.  Powers' definition differs somewhat from 
Definition 1.5, in that he requires only weak convergence 
to $\omega$
$$
\lim_{t\to\infty}\rho(\alpha_t(x))=\omega(x), \qquad x\in\Cal B(H), 
$$
for every normal state $\rho$.  But as the following observation 
shows, the two definitions are in fact equivalent.  
\endremark

\proclaim{Proposition 1.6}
Let $\{\rho_i: i\in I\}$ be a net of normal states of $M=\Cal B(H)$
and let $\omega$ be a normal state such that 
$$
\lim_i\rho_i(x)=\omega(x), \tag{1.6.1}
$$
for every compact operator $x$.  Then $\lim_i\|\rho_i-\omega\| = 0$.  
\endproclaim
\demo{proof}
Choose $\epsilon>0$.  Since $\omega$ is a normal s we can find
a finite rank projection $p$ such that 
$$
\omega(p)\geq 1-\epsilon.  \tag{1.7}
$$
Since $pMp\cong \Cal B(pH)$ is a finite dimensional space of 
finite-rank operators, (1.6.1) implies that we have norm convergence
$$
\lim_i\|\rho_i\restriction_{pMp} -\omega\restriction_{pMp}\|=0,
$$
and hence
$$
\sup_{x\in M, \|x\|\leq 1}|\rho_i(pxp)-\omega(pxp)| \to 0, \tag{1.8}
$$
as $i\to\infty$.  Now in general, we have 
$$
\|\rho_i-\omega\|\leq \sup_{\|x\|\leq 1}|\rho_i(pxp)-\omega(pxp)| +
\sup_{\|x\|\leq 1}|\rho_i(x-pxp)| + \sup_{\|x\|\leq 1}|\omega(x-pxp)|.
$$
By (1.8), the first term on the right tends to $0$ as $i\to\infty$, 
and we can estimate the second and third terms as follows.  
Writing $x-pxp=(\bold 1-p)x+px(\bold 1-p)$, we find from the 
Schwarz inequality that 
$$
|\rho_i((\bold 1-p)x)|^2\leq \rho_i(\bold 1-p)\rho_i(x^*x)\leq
(1-\rho_i(p))\|x\|^2
$$
and hence 
$$
|\rho_i((\bold 1-p)x)|\leq (1-\rho_i(p))^{1/2}\|x\|.  
$$
Similarly, 
$$
|\rho_i(px(\bold 1-p))|\leq (1-\rho_i(p))^{1/2}\|px\|\leq 
(1-\rho_i(p))^{1/2}\|x\|.  
$$
It follows that 
$$
\sup_{\|x\|\leq 1}|\rho_i(x-pxp)|\leq 2(1-\rho_i(p))^{1/2}.  
$$
Since $1-\rho_i(p)$ tends to $1-\omega(p)\leq \epsilon$ as
$i\to\infty$, it follows that
$$
\limsup_{i\to\infty}\sup_{\|x\|\leq 1}|\rho_i(x-pxp)|\leq 
2 \epsilon^{1/2}.  
$$
Similar estimates show that 
$$
\sup_{\|x\|\leq 1}|\omega(x-pxp)|\leq 2 \epsilon^{1/2}.  
$$
using (1.8), we conclude that 
$$
\limsup_{i\to\infty}\|\rho_i-\omega\|\leq 4 \epsilon^{1/2}
$$
and (1.6.1) follows because $\epsilon$ is arbitrary \qed
\enddemo

\remark{Remarks}
Suppose that $\alpha=\{\alpha_t: t\geq 0\}$ is a pure 
\esg\ acting on an arbitrary von Neumann algebra $M$, and 
that $\omega$ is a normal state of $M$ which is invariant 
under $\alpha$.  Then for every normal state $\rho$, 
Proposition 1.1 implies that 
$$
\lim_{t\to\infty}\|\rho\circ\alpha_t-\omega\| = 
\lim_{t\to\infty}\|\rho\circ\alpha_t-\omega\circ\alpha_t\| = 0, 
$$
hence $\omega$ is an absorbing state.  Conversely, if an 
\esg\ $\alpha$ has an absorbing state, then by Proposition
1.1 $\alpha$ must 
be a pure \esg.  Thus we have the following description of the
relationship between absorbing states and pure \esg s: 
\endremark

\proclaim{Proposition 1.9}
Let $\alpha=\{\alpha_t: t\geq 0\}$ be an \esg\ acting on a 
von Neumann algebra $M$ which has a normal 
invariant state $\omega$.  Then $\alpha$ is pure 
if and only if $\omega$ is an absorbing state.  
\endproclaim

\remark{Remarks}
Every abelian semigroup is amenable.  Thus one can make use 
of a Banach limit on the additive semigroup 
of nonnegative reals to average any \esg\ in the point-weak 
operator topology to show that there is a state of 
$\Cal B(H)$ which is invariant under the action of 
the \esg.  However, invariant states constructed by such 
devices tend to be singular.  Indeed, 
the results of \cite{6} show that there are pure \esg s 
(acting on $\Cal B(H)$) which do not have normal invariant states.  
\endremark

%\pagebreak

\subheading{2.  Pure CP semigroups}
\proclaim{Definition 2.1}
A CP semigroup is a semigroup $P=\{P_t: t\geq 0\}$ of 
normal completely positive maps of $\Cal B(H)$ which 
satisfies the natural continuity property
$$
\lim_{t\to 0+}\<P_t(x)\xi,\eta\>=\<x\xi,\eta\>, 
\qquad x\in\Cal B(H), \xi,\eta\in H.  
$$
$P$ is called unital if $P_t(\bold 1)=\bold 1$ for 
every $t\geq 0$.  

A unital CP semigroup $P$ is said to be pure if, for 
every pair of normal states $\rho_1, \rho_2$ of 
$\Cal B(H)$ we have 
$$
\lim_{t\to\infty}\|\rho_1\circ P_t-\rho_2\circ P_t\| = 0.  
$$
\endproclaim

Notice that pure CP semigroups are required to be unital.  
Unital CP semigroups are often called {\it quantum 
dynamical semigroups} in the mathematical physics 
literature.  The purpose of this section is to briefly 
discuss the relationship between pure CP semigroups 
and pure \esg s.  This relationship is not bijective, but it is 
close enough to being so that results in one category
usually have immediate implications for the other.  

For example, suppose that $P$ is a pure CP semigroup 
acting on $\Cal B(H)$.  A recent dilation theorem of 
B. V. R. Bhat \cite{7,8} implies that there is a Hilbert
space $K\supseteq H$ and an \esg\ 
$\alpha=\{\alpha_t: t\geq 0\}$ acting on $\Cal B(K)$ which 
is a {\it dilation} of $P$ in the following sense.  Letting 
$p_0\in\Cal B(K)$ be the projection onto the subspace $H$ 
and identifying $\Cal B(H)$ with the hereditary subalgebra
$M_0=p_0\Cal B(K)p_0$ of $\Cal B(K)$,  
then we have 
$$
\align
\alpha_t(p_0)&\geq p_0, \qquad\qquad \text{and}\tag{2.2.1}\\
P_t(x) &= p_0\alpha_t(x)p_0, \qquad x\in M_0 \tag{2.2.2}
\endalign
$$
for every $t\geq 0$.  Because of (2.2.1), the operator 
$$
p_\infty = \lim_{t\to\infty}\alpha_t(p_0)
$$
exists as a strong limit of projections, and is therefore 
a projection fixed under the action of $\alpha$.  By compressing
$\alpha$ to the hereditary subalgebra $p_\infty\Cal B(K)p_\infty$ 
if necessary, we can assume that $K=p_\infty K$ and hence 
that 
$$
\alpha_t(p_0)\uparrow\bold 1_K, \qquad\text{as }t\to\infty.  \tag{2.3}
$$
When (2.3) is satisfied we will say that $\alpha$ is a {\it dilation}
of $P$.  

Dilations in this sense are not unique.  In order to obtain 
uniqueness (up to conjugacy), one must in general compress 
$\alpha$ to a smaller hereditary subalgebra of $\Cal B(K)$.  
Once that is done  
$\alpha$ is called a {\it minimal} dilation of $P$.  The 
issue of minimality is a subtle one, and we will not have to 
be very specific about its nature here (see \cite{2} for more
detail).  For our purposes, it is enough to know that every 
dilation can be compressed uniquely to a minimal dilation, and 
that minimal dilations are unique up to conjugacy.  
Moreover, nonminimal dilations of a given CP semigroup 
exist in profusion.  For example, the trivial CP semigroup 
acting on $\Bbb C$ has many dilations to nontrivial, 
nonconjugate \esg s \cite{13}.  The 
following result implies that all such dilations are pure.  

\proclaim{Proposition 2.4}
Let $P=\{P_t: t\geq 0\}$ be a pure CP semigroup
acting on $\Cal B(H)$.  
Then every dilation of $P$ to an \esg\ is pure.  
\endproclaim
\demo{proof}
Let $\alpha$ be a dilation of $P$ which acts on $\Cal B(K)$,
$K$ being a Hilbert space containing $H$.   Letting $p_0\in\Cal B(K)$
be the projection on $H$, then by (2.3) we see that the subspaces
$$
K_t = \alpha_t(p_0)K
$$
increase with $t$ and their union is dense in $K$.  If we let 
$\Cal N_t$ denote the set of all normal states $\rho$ of $\Cal B(K)$ 
which can be represented in the form 
$$
\rho(x) = \sum_k\<x\xi_k,\xi_k\>, 
$$
with vectors $\xi_1, \xi_2, \dots \in K_t$, then the sets 
$\Cal N_t$ increase with $t$ and their union is {\it norm}-dense 
in the space of all normal states of $\Cal B(K)$.  

Using this observation together with Proposition 1.1, it is 
enough to show that for every $t>0$ and every pair of normal 
states $\rho_1$, $\rho_2\in \Cal N_t$, we have 
$$
\lim_{s\to\infty}\|\rho_1\circ\alpha_s-\rho_2\circ\alpha_s\|=0.  \tag{2.5}
$$
To prove (2.5), fix $t>0$ and choose $s>t$.  We claim that 
for $k=1,2$ and $x\in \Cal B(K)$ we have 
$$
\rho_k(\alpha_s(x)) = \rho_k(\alpha_t(P_{s-t}(p_0xp_0))).  \tag{2.6}
$$
Indeed, since $p_0\leq \alpha_{s-t}(p_0)$ we have 
$$
P_{s-t}(p_0xp_0)=p_0\alpha_{s-t}(p_0xp_0) =
p_0\alpha_{s-t}(p_0)\alpha_{s-t}(x)\alpha_{s-t}(p_0)p_0 = 
p_0\alpha_{s-t}(x)p_0,   
$$
so that  
$$
\alpha_t(P_{s-t}(p_0xp_0)))=\alpha_t(p_0\alpha_{s-t}(x)p_0) =
\alpha_t(p_0)\alpha_s(x)\alpha_t(p_0).  
$$
Hence the right side of (2.6) can be written 
$$
\rho_k(\alpha_t(p_0)\alpha_s(x)\alpha_t(p_0)).  
$$
Since $\rho_k$ belongs to $\Cal N_t$ we must have 
$\rho_k(\alpha_t(p_0)z\alpha_t(p_0))=\rho_k(z)$ for every 
$z\in\Cal B(K)$, and (2.6) follows.  

Letting $\sigma_k$ be the restriction of $\rho_k\circ\alpha_t$ to 
$M_0=p_0\Cal B(K)p_0$ we find that for every $x\in\Cal B(K)$, 
 
$$
|\rho_1(\alpha_s(x))-\rho_2(\alpha_s(x))|=
|\sigma_1(P_{s-t}(p_0xp_0)-\sigma_2(P_{s-t}(p_0xp_0)|.  
$$
Thus 
$$
\|\rho_1\circ\alpha_s-\rho_2\circ\alpha_s\| =
\|\sigma_1\circ P_{s-t}-\sigma_2\circ P_{s-t}\|  
$$
must tend to $0$ as $s$ 
tends to $\infty$, and (2.5) follows\qed
\enddemo

Suppose now that we start with a pure \esg\ acting on 
$\Cal B(H)$.  It is not always possible to locate a 
CP semigroup as a compression of $\alpha$ because we know 
of no general method for locating a projection $p_0\in\Cal B(H)$
satisfying $\alpha_t(p_0)\geq p_0$ for every $t$.  However, if 
$\alpha$ has an invariant normal state $\omega$, 
then the support projection of $\omega$ provides such a 
projection $p_0$.  To see that, simply notice that 
$\omega\circ\alpha_t(\bold 1-p_0)=\omega(\bold 1-p_0)=0$, 
hence $\alpha_t(\bold 1-p_0)\leq \bold 1-p_0$,
hence $\alpha_t(p_0)\geq p_0$.  

Given such a projection $p_0$, we can compress $\alpha$ to 
obtain a family of normal completely positive maps 
$P=\{P_t: t\geq 0\}$ of $\Cal B(p_0H)\cong p_0\Cal B(H)p_0$
by way of 
$$
P_t(x) = p_0\alpha_t(x)p_0, \qquad t\geq 0, x\in p_0\Cal B(H)p_0.  
\tag{2.7} 
$$
The fact that $\alpha_t(p_0)\geq p_0$ insures that $P$ is 
in fact a CP semigroup.  The following summarizes these 
remarks.  

\proclaim{Proposition 2.8}
Suppose that $\alpha$ is a pure \esg\ acting on $\Cal B(H)$ and
$\omega$ is a normal $\alpha$-invariant state with support projection 
$p_0$.  Then the CP semigroup P defined by (2.7) is pure, and the 
restriction $\omega_0$ of $\omega$ to $p_0\Cal B(H)p_0\cong \Cal B(p_0H)$
is a faithful normal $P$-invariant state which is absorbing in 
the sense that for every normal state $\rho$ of $\Cal B(p_0H)$, 
$$
\lim_{t\to\infty}\|\rho\circ P_t-\omega_0\| = 0.  
$$

If $\omega$ is weakly continuous and not a pure state of $\Cal B(H)$, 
then $P$ may be considered a CP semigroup acting on a matrix 
algebra $M_n(\Bbb C)$, $n=2,3,\dots$.  
\endproclaim

The preceding discussion shows the extent to which the theory of 
pure \esg s having an absorbing state can be reduced to the theory
of CP semigroups having a {\it faithful} absorbing state.  
While the latter problem
is an attractive one in general, we still 
lack tools that are appropriate
for arbitrary invariant normal states.  The following sections
address the case of weakly continuous invariant states.

\subheading{3.  Perturbations and invariant states}

In order to describe the pure CP semigroups acting on matrix
algebras we must first obtain information about invariant 
states.  More precisely, given a {\it faithful} state 
$\omega$ on a matrix algebra $M=M_N(\Bbb C)$, 
$N=2,3,\dots$, we want to identify the unital CP semigroups
$P=\{P_t: t\geq 0\}$ that leave $\omega$ invariant in 
the sense that 
$$
\omega\circ P_t =\omega, \qquad t\geq 0.  
$$
It is not obvious that such semigroups exist when 
$\omega$ is not a tracial state.  In this section
we characterize the generators of such semigroups 
up to perturbations (Theorem 3.8) and we 
give explicit examples in Corollary 3.16.  

In general, the generator $L$ of a CP semigroup has
a decomposition of the form 
$$
L(x) = P(x) + kx + xk^*, \qquad x\in M \tag{3.1}
$$
where $P$ is a completely positive map on $M$ and 
$k\in M$ \cite{10}.  The associated semigroup 
$\{\exp{tL}: t\geq 0\}$ is unital iff 
$$
L(\bold 1) = 0   \tag{3.2}
$$
and it leaves $\omega$ invariant iff
$$
\omega\circ L = 0.  \tag{3.3}
$$
It is easy to satisfy (3.2), but less easy to 
satisfy both (3.2) and (3.3).  Indeed, setting $x=\bold 1$
in (3.1) we find that (3.2) holds iff $k$ has a 
Cartesian decomposition 
$$
k=-1/2P(\bold 1) + \ell,
$$
where $\ell$ is an element of $M$ satisfying $\ell^*=-\ell$. 
In this case (3.1) becomes 
$$
L(x) = P(x) - 1/2(P(\bold 1)x+xP(\bold 1)) + [\ell,x].  \tag{3.4}
$$

There is a natural decomposition of this operator 
corresponding to the Cartesian decomposition of $k$:
$$
L(x) = L_0(x) +[\ell,x],
$$
where $L_0$ is the ``unperturbed" part of $L$
$$
L_0(x) = P(x) - 1/2(P(\bold 1)x + xP(\bold 1)).  \tag{3.5}
$$
Notice that both $L_0$ and $L$ generate unital CP 
semigroups, and because of (3.3) the semigroup generated
by $L$ leaves $\omega$ invariant.  If $\omega$ is not a 
trace then the unperturbed CP semigroup 
$\exp{tL_0}$ need not leave $\omega$ invariant
(see Proposition 3.18).  
Thus we are led to seek {\it perturbations}
of $L_0$ which solve both equations (3.2) and (3.3).  

In order to discuss this issue in more concrete terms, 
let $\Omega$ be the density matrix of the state $\omega$, 
$$
\omega(x) = \tr (\Omega x), \qquad x\in M.  
$$
Since $\omega$ is faithful, $\Omega$ is a positive 
invertible operator.  More generally, we identify 
the dual $M^\prime$ of $M$ with $M$ itself in the usual way,
the isomorphism $a\in M\mapsto \omega_a\in M^\prime$ 
being defined by 
$$
\omega_a(x) = \tr (ax), \qquad x\in M.  
$$
For every linear map $L: M\to M$ the dual 
map $L_*$, defined on $M^\prime$ by
$L_*(\rho) = \rho\circ L$, becomes 
$$
\tr(L_*(y)x) = \tr(yL(x)), \qquad x,y\in M.  
$$

Now a linear map $L: M\to M$ satisfies 
$\omega\circ L=0$ iff its dual satisfies 
$L_*(\Omega) = 0$.  If we choose a completely positive
map $P: M\to M$ and define $L_0$ as in (3.5), then we 
seek a skew-adjoint operator $\ell\in M$ satisfying 
the operator equation 
$$
L_{0*}(\Omega) = \ell \Omega-\Omega \ell.  \tag{3.6}
$$
It is not always possible to solve (3.6).  But if a 
solution $\ell_0$ exists then there are infinitely many, the 
most general one having the form $\ell = \ell_0 + k$, 
$k$ being a skew-adjoint operator commuting with $\Omega$.  

We will show that (3.6) is solvable iff 
$P$ satisfies a certain symmetry requirement.  
The symmetry involves an involution $\#$ and is 
described as follows.  For every linear map $L: M\to M$, 
let $L^\#: M\to M$ be the linear map 
$$
L^\#(x) = \Omega^{-1/2}L_*(\Omega^{1/2}x\Omega^{1/2})\Omega^{-1/2}.
\tag{3.7}
$$
For our purposes, the important properties of the 
operation $L\mapsto L^\#$ are summarized as follows.  

\proclaim{Proposition}
$L\mapsto L^\#$ is a linear isomorphism satisfying 
$L^{\#\#}=L$, and if $L$ is completely positive then
so is $L^\#$.  
\endproclaim
\demo{sketch of proof}
The argument is completely straightforward.  A direct 
computation shows that 
$$
(L^\#)_*(x) = \Omega^{1/2}L(\Omega^{-1/2}x\Omega^{-1/2})\Omega^{1/2},
$$
from which $L^{\#\#}=L$ is immediate.  The fact that 
$\#$ preserves complete positivity follows from the 
fact that if $P$ is a completely positive map then 
so is $P_*$.\qed
\enddemo

\proclaim{Theorem 3.8}
Let $\omega$ be a faithful state on a matrix algebra
$M$, let $Q: M\to M$ be a completely positive linear 
map, and define $Q^\#$ by (3.7).  Then the following 
are equivalent.  
\roster
\item"{(i)}"
There is a unital CP semigroup $P=\{P_t: t\geq 0\}$ which 
leaves $\omega$ invariant and whose generator has the form
$$
L(x) = Q(x) + kx + xk^*
$$
for some $k\in M$.  
\item"{(ii)}"
For every minimal spectral projection $e$ of $\Omega$ 
we have $eQ(\bold 1)e=eQ^\#(\bold 1)e$.  
\endroster
\endproclaim

Our proof of Theorem 3.8 is based on the following 
general result.  Let $A$ be the centralizer algebra of 
$\omega$, 
$$
A = \{a\in M: \omega(ax) = \omega(xa), x\in M\}.  
$$
If we consider the spectral decomposition of $\Omega$, 
$$
\Omega = \sum_{k=1}^r\lambda_k e_k
$$
where $e_1, \dots, e_r$ are the minimal spectral projections
of $\Omega$ and $0<\lambda_1<\dots <\lambda_r$ are the 
distinct eigenvalues, then $A$ is the commutant of 
$\{\Omega\}$ and hence 
$$
A=\{a\in M: ae_k=e_ka, 1\leq k\leq r\}.  
$$
$A$ is a direct sum of full matrix algebras, and the 
restriction of $\omega$ to $A$ is a faithful tracial 
state.  The natural conditional expectation 
$E_A: M\to A$ is given by
$$
E_A(x) = \sum_ke_kxe_k, \qquad x\in M.  
$$
The following result implies that the solvability of 
equation (3.6) depends only on the compression of 
$L$ to the centralizer algebra $A$.  

\proclaim{Lemma 3.9}
Let $\omega$ be a faithful state of $M$ and let 
$L: M\to M$ be a linear map satisfying 
$L(x)^*=L(x^*)$, $x\in M$.  The following are 
equivalent:
\roster
\item"{(i)}"
There is a skew-adjoint operator $\ell\in M$ such 
that the perturbation $L^\prime(x) = L(x)+[\ell,x]$
satisfies $\omega\circ L^\prime = 0$.  
\item"{(ii)}"
The restriction of $\omega\circ L$ to $A$ vanishes.  
\endroster
More generally, setting $L_0=E_ALE_A$, there is 
a perturbation 
$L^\prime(x)=L(x)+[\ell,x]$ of the form (i) such
that $\omega\circ L^\prime = \omega\circ L_0$.  
\endproclaim
\demo{proof: (i)$\implies$(ii)}  Suppose that $\ell$ is 
an operator in $M$ for which $\omega\circ L^\prime=0$, 
$L^\prime$ being the operator of part (i).  Since 
$\omega(\ell a-a\ell)=0$ for all $a$ in the centralizer
algebra we have 
$$
\omega(L(a))=\omega(L(a)+[\ell,a])=\omega\circ L^\prime(a)=0,
$$ 
hence (ii).  

We now prove the general assertion of the last sentence.  
Noting that $\omega\circ E_A=\omega$, we have 
$$
\omega(L_0(x))=\omega(L(E_A(x))), \qquad x\in M
$$ 
and hence we must exhibit an operator 
$\ell\in M$ satisfying $\ell^*=-\ell$ and 
$$
\omega(L(x)+[\ell,x] -L(E_A(x)))=0, \qquad x\in M.  
$$
After dualizing, the previous equation becomes 
$$
L_*(\Omega)-[\ell,\Omega] - E_A(L_*(\Omega))=0,
$$
or 
$$
L_*(\Omega)-E_A(L_*(\Omega))=\ell\Omega-\Omega\ell.
\tag{3.10}
$$
Let $T$ be the left side of (3.10).  $T$ is a self-adjoint
operator satisfying $E_A(T)=0$.  Thus if 
$$
\Omega = \sum_{k=1}^r\lambda_ke_k
$$
is the spectral decomposition of $\Omega$ then we have 
$e_kTe_k=0$ for all $k$.  Set 
$$
\ell = \sum_{i\neq j}\frac{1}{\lambda_j-\lambda_i}e_iTe_j.  
$$
It is obvious that $\ell^*=-\ell$, and since 
$\Omega e_k=e_k\Omega=\lambda_ke_k$ for all $k$ we have 
$$
\align
\Omega\ell &= 
\sum_{i\neq j}\frac{\lambda_i}{\lambda_j-\lambda_i}e_iTe_j, \\
\ell\Omega &= 
\sum_{i\neq j}\frac{\lambda_j}{\lambda_j-\lambda_i}e_iTe_j.  
\endalign
$$
Hence 
$$
\ell\Omega-\Omega\ell = \sum_{i\neq j}e_iTe_j = T,
$$
as required.  

The implication (ii)$\implies$(i) follows immediately, 
for if $\omega\circ L(a)=0$ for all $a\in A$, then because
$\omega\circ E_A = \omega$ we have 
$\omega\circ L_0=0$.  Thus the preceding argument 
gives a perturbation $L^\prime$ of the form (i) satisfying 
$\omega\circ L^\prime = \omega\circ L_0 =0$\qed
\enddemo

\demo{proof of Theorem 3.8}
Let $Q$ be a completely positive map
and define $L: M\to M$ by 
$$
L(x) = Q(x)-1/2(Q(\bold 1)x+xQ(\bold 1)).  
$$
The assertion (i) of Theorem 3.8 is equivalent to 
the existence of a  skew-adjoint operator 
$\ell\in M$ such that 
$$
\omega(L(x)+[\ell,x])=0, \qquad x\in M.  \tag{3.11}
$$
By Lemma 3.9, the latter is equivalent to 
$$
\omega(L(a))=0, \qquad a\in A.  \tag{3.12}
$$
Thus we have to show that (3.12) is equivalent to the 
operator equation 
$$
E_A(Q(\bold 1))=E_A(Q^\#(\bold 1)).  \tag{3.13}
$$
Looking first at (3.12), we have 
$$
\omega(L(a))=\omega(Q(a))-1/2\omega(Q(\bold 1)a+aQ(\bold 1)).  
$$
Now since every element $a\in A$ commutes with $\Omega$ we have
$$
\align
1/2\omega(Q(\bold 1)a+aQ(\bold 1))&=
1/2\tr(\Omega Q(\bold 1)a+\Omega aQ(\bold 1))\\
&=
\tr(\Omega Q(\bold 1)a) = \omega(Q(\bold 1)a).
\endalign
$$
Hence (3.12) asserts that 
$$
\omega(Q(a))-\omega(Q(\bold 1)a) = 0, \qquad a\in A.  \tag{3.14}
$$
Writing 
$$
\align
\omega(Q(a))&=\tr(\Omega Q(a))=\tr(Q_*(\Omega)a) \\
&= 
\tr(\Omega\cdot\Omega^{-1/2}Q_*(\Omega^{1/2}\cdot\Omega^{1/2})\Omega^{-1/2}a)
=\omega(Q^\#(\bold 1)a),
\endalign
$$
we rewrite (3.14) as 
$$
\omega((Q^\#(\bold 1)-Q(\bold 1))a)= 0, \qquad a\in A.  
$$
Since $\omega\circ E_A=\omega$ and $E_A(xa)=E_A(x)a$ for 
$a\in A$ the preceding formula becomes
$$
\omega(E_A(Q^\#(\bold 1)-Q(\bold 1))a)=0,\qquad a\in A.  
$$
Since $\omega\restriction_A$ is a faithful trace on $A$, 
the latter is equivalent to equation (3.13).\qed
\enddemo

\remark{Remark 3.15}
In the important case where $\omega$ is the tracial state
on $M$ the density matrix of $\omega$ is a scalar,  
the map $\#$ reduces to the dual mapping $L^\# = L_*$, 
and $E_A$ is the identity map.  In this case the criterion 
(ii) of Theorem 3.8 degenerates to $Q(\bold 1) = Q_*(\bold 1)$.  
For example, if $Q$ has the form 
$$
Q(x)=\sum_{k=1}^rv_kxv_k^*
$$
where $v_1,v_2,\dots,v_r\in M$, then condition (ii) 
becomes 
$$
\sum_{k=1}^rv_kv_k^*=\sum_{k=1}^rv_k^*v_k.  
$$
Moreover, when this condition is satisfied and 
$\omega$ is the tracial state no perturbation is
necessary.  One simply shows by a direct calculation
that the mapping 
$$
L(x) = Q(x) -1/2(Q(\bold 1)x+xQ(\bold 1))
$$
satisfies $\tr\circ L=0$ iff $Q(\bold 1) = Q_*(\bold 1)$.  
\endremark

\proclaim{Corollary 3.16}
Let $\omega$ be a faithful state on $M$ with density 
matrix $\Omega$ and let $v_1,\dots,v_r\in M$ satisfy
$$
\sum_{k=1}^rv_kv_k^*=\sum_{k=1}^rv_k^*v_k.  
$$
Then there is a unital $\omega$-preserving CP semigroup
whose generator has the form 
$$
L(x) = \Omega^{-1/2}(\sum_{k=1}^rv_kxv_k^*)\Omega^{-1/2}
+kx + xk^*
$$
for some operator $k\in M$.  
\endproclaim
\demo{proof}
Let $Q$ be the completely positive map
$$
Q(x)=\Omega^{-1/2}(\sum_{k=1}^rv_kxv_k^*)\Omega^{-1/2}.  
$$
By Theorem 3.8 it suffices to show that 
$Q^\#(\bold 1)=Q(\bold 1)$.  A direct computation shows
that the dual of $Q$ is given by
$$
Q_*(x) = \sum_{k-1}^rv_k^*\Omega^{-1/2}x\Omega^{-1/2}v_k
$$
Hence 
$$
Q^\#(\bold 1) = \Omega^{-1/2}Q_*(\Omega)\Omega^{-1/2} =
\Omega^{-1/2}(\sum_{k=1}^rv_k^*v_k)\Omega^{-1/2}.  
$$
The right side is $Q(\bold 1)$ because of the hyposthesis
on $v_1,\dots,v_r$.\qed
\enddemo

\remark{Remark 3.17: The necessity of perturbations}
In view of Remark 3.15 it is natural to ask if nontrivial  
perturbations are really necessary, and we conclude this
section with some remarks concerning that issue.  
Suppose that $P$ is a normal completely positive 
map of $M$ and $L$ is the unperturbed generator
$$
L(x)=P(x)-1/2(P(\bold 1)x+xP(\bold 1)).  \tag{3.18}
$$
\endremark

\proclaim{Proposition 3.19}
Let $\omega$ be a faithful state on $M=M_N(\Bbb C)$ which 
is not a trace.  Then there is an operator $L$ of the 
form (3.18) and a skew-adjoint operator $\ell\in M$ such 
that if $L^\prime(x)=L(x)+[\ell,x]$ then $\omega\circ L\neq 0$
while $\omega\circ L^\prime = 0$.  
\endproclaim

\demo{proof}
Consider the spectral decomposition of the density 
matrix of $\omega$
$$
\Omega = \sum_{k=1}^r\lambda_ke_k.
$$
We must have $r\geq 2$ because $\omega$ is not a trace.  
Choose a nonzero partial isometry $v$ satisfying 
$v^*v\leq e_1$ and $vv^*\leq e_2$.  Since $\Omega$ is 
an invertible positive operator there is an $\epsilon>0$ 
such that 
$$
\Omega^\prime = \Omega + \epsilon(v+v^*)
$$
is positive.  Since the trace of $\Omega^\prime$ is $1$
we may consider the state $\omega^\prime$ having density
matrix $\Omega^\prime$.  Let $P$ be a normal 
completely positive map satisfying $P(\bold 1)=\bold 1$
and $\omega\circ P=\omega^\prime$ (there are many such maps, 
the simplest one being $P(x) = \omega^\prime(x)\bold 1$),  
and define 
$$
L(x)=P(x)-x.  
$$
Then $\omega\circ L=\omega^\prime-\omega\neq 0$.  

On the other hand, since 
$P_*(\Omega)=\Omega^\prime$ we have 
$$
P^\#(\bold 1)=\Omega^{-1/2}P_*(\Omega)\Omega^{-1/2}=
\Omega^{-1/2}\Omega^\prime\Omega^{-1/2}.  
$$
Thus, letting 
$$
E_A(x) = \sum_{k=1}^r e_kxe_k
$$
be the conditional expectation onto the centralizer algebra
of $\omega$ and using $E_A(v)=E_A(v^*)=0$, we have 
$E_A(\Omega^\prime)=E_A(\Omega)$.  Hence  
$$
E_A(P^\#(\bold 1))=
\Omega^{-1/2}E_A(\Omega^\prime)\Omega^{-1/2} = \bold 1.  
$$
From Theorem 3.8 we may conclude
that there is a skew-adjoint operator 
$\ell$ such that the perturbation 
$$
L^\prime(x) = L(x)+[\ell,x]
$$
satisfies $\omega\circ L^\prime=0$\qed
\enddemo

\subheading{4.  Ergodicity and purity}
The purpose of this section is to give a concrete 
characterization of the generators of pure CP semigroups
acting on matrix algebras, given that the CP semigroup 
has a faithful invariant state (Theorem 4.4).  
\proclaim{Definition 4.1}
A unital CP semigroup $P=\{P_t: t\geq 0\}$ acting on 
$\Cal B(H)$ is called ergodic if the only operators
$x$ satisfying $P_t(x)=x$ for every $t\geq 0$ are 
scalars.  
\endproclaim

The set $\Cal A=\{x\in\Cal B(H): P_t(x)=x, t\geq 0\}$ 
is obviously a weak$^*$-closed self-adjoint linear 
subspace of $\Cal B(H)$ containing the identity.  
In general it need not be a von Neumann algebra, but 
as we will see presently, it is a von Neumann 
algebra in the cases of primary interest for our
purposes here.  

\proclaim{Proposition 4.2}
Every pure CP semigroup is ergodic.  
\endproclaim 
\demo{proof}
Suppose $P=\{P_t: t\geq 0\}$ is pure and $x$ is 
an operator satisfying $\|x\|\leq 1$ and 
$P_t(x)=x$ for every 
$t$.  To show that $x$ must be a scalar multiple
of $\bold 1$ it suffices to show that for every 
normal linear functional $\rho$ on $\Cal B(H)$ 
satisfying $\rho(\bold 1)=0$ we have 
$\rho(x)=0$.  Since any normal linear functional
$\rho$ satisfying $\rho(\bold 1)=0$ can be decomposed
into a sum of the form 
$$
\rho = b(\rho_1-\rho_2) +ic(\rho_3-\rho_4)
$$
where $b$ and $c$ are real numbers and the 
$\rho_k$ are normal states, we conclude from 
the purity of $P$ that 
$$
\lim_{t\to\infty}\|\rho\circ P_t\|=0.  
$$
Since $x$ is fixed under the action of $P$ 
we have 
$$
|\rho(x)|=|\rho(P_t(x))|\leq \|\rho\circ P_t\|
$$
for every $t\geq 0$, from which 
$\rho(x)=0$ follows.\qed
\enddemo

\proclaim{Proposition 4.3}
Let $P=\{P_t: t\geq 0\}$ be a unital CP 
semigroup which leaves invariant some faithful 
normal state of $\Cal B(H)$.  Then 
$$
\Cal A=\{a\in \Cal B(H): P_t(a)=a, t\geq 0\}
$$ is a von Neumann algebra.  
Assuming further that 
$P$ has a bounded generator $L$ 
represented in the form 
$$
L(x) = \sum_j v_jxv_j^* + kx + xk^*  \tag{4.3.1}
$$
for operators $k,v_1,v_2,\dots\in \Cal B(H)$, then
$\Cal A$ is the commutant of the von Neumann algebra
generated by $\{k,v_1, v_2, \dots\}$.  
\endproclaim

\demo{proof}
In view of the preceding remarks, the first paragraph
will follow if we show that $\Cal A$ is closed under 
operator multiplication.  By polarization, it is 
enough to show that $a\in\Cal A\implies a^*a\in\Cal A$.  
For each $a\in\Cal A$ we have by the Schwarz inequality
$$
a^*a=P_t(a)^*P_t(a)\leq P_t(a^*a)
$$
for every $t\geq 0$.  Letting $\omega$ be a faithful 
state invariant under $P$ we have 
$\omega(P_t(a^*a)-a^*a)=0$, and hence $P_t(a^*a)=a^*a$.  
Thus $a^*a\in\Cal A$.  

Suppose now that $P$ has a bounded generator of the 
form (4.3.1), and let $\Cal B$ be the $*$-algebra 
generated by $\{k,v_1,v_2,\dots\}$.  Noting that
$\Cal A=\{x\in M: L(x)=0\}$, we show 
that $\Cal A=\Cal B^\prime$.  
If $x\in\Cal B^\prime$ then (4.3.1) becomes
$$
L(x)=x(\sum_jv_jv_j^* + k + k^*)=xL(\bold 1)=0.  
$$
It follows that $\exp{tL}(x)=x$ for every 
$t$, hence $x\in\Cal A$.  

For the inclusion $\Cal A\subseteq\Cal B^\prime$, 
we claim first that for every $a\in \Cal A$, 
$$
[v_j,a]=v_ja-av_j=0, \qquad j=1,2,\dots.  
$$
Indeed, since $\bold 1,a,a^*$, and $aa^*$ all belong
to $\Cal A$ and $L(\Cal A)=\{0\}$, we have 
$$
L(aa^*)-aL(a^*)-L(a)a^*+aL(\bold 1)a^*=0.  
$$
Substituting the formula (4.3.1) for $L$ in the 
above we find that the terms involving $k$ drop 
out and we are left with the formula 
$$
\sum_k[v_j,a][v_j,a]^* =-\sum_j[v_j,a][v_j^*,a^*]=0.  
$$
It follows that $[v_j,a]=0$ for every $k$.  Replacing 
$a$ with $a^*$ we see that $a$ must commute with the 
self-adjoint set of operators 
$\{v_1,v_2,\dots,v_1^*,v_2^*,\dots\}$.  

Now since $L(\bold 1)=0$, it follows from (4.3.1) 
that $\sum_jv_jv_j^*+k+k^*=0$, and hence 
$k$ has Cartesian decomposition 
$k=-h+\ell$ where 
$$
h=1/2\sum_jv_jv_j^*
$$
and $\ell$ is a skew-adjoint operator.  Setting 
$$
L_0(x) = \sum_jv_jxv_j^* -hx - xh, 
$$
we have 
$$
L(x)=L_0(x)+[\ell,x], 
$$
and $L_0(\Cal A)=\{0\}$ by what was just 
proved.  Thus, for $a\in\Cal A$ 
$$
[\ell,a]=L(a)=0,
$$
and hence $a$ must commute with $\ell$ as well.  The
inclusion $\Cal A\subseteq \Cal B^\prime$ follows.\qed
\enddemo

\proclaim{Theorem 4.4}
Let $P=\{P_t: t\geq 0\}$ be a unital 
CP semigroup acting on a matrix algebra 
$M=M_N(\Bbb C)$, $N=2,3,\dots$ which leaves 
invariant some faithful state $\omega$.  Let
$$
L(x)=\sum_{j=1}^rv_jxv_j^*+kx+xk^*
$$
be the generator of $P$.  Then the following 
are equivalent:
\roster
\item"{(i)}"
$P$ is pure.  
\item"{(ii)}"
$P$ is ergodic.  
\item"{(iii)}"
The set of operators 
$\{k,k^*,v_1,\dots,v_r,v_1^*,\dots,v_r^*\}$
is irreducible.  
\endroster
\endproclaim

\demo{proof}
In view of Propositions 4.2 and 4.3, we need only 
prove the implication (ii)$\implies$(i).  Assuming 
that $P$ is ergodic, we consider its generator 
$L$ as an operator on the Hilbert space 
$L^2(M,\omega)$ with inner product
$$
\<x,y\>=\omega(y^*x), \qquad x,y\in M.  
$$
We have $L(\bold 1)=0$ because $P$ is unital, 
and $L^*(\bold 1)=0$ follows from the fact that
$\omega\circ L=0$, $L^*$ denoting the adjoint 
of $L\in\Cal B(L^2(M,\omega))$.  It follows that
$\{\lambda\bold 1: \lambda\in\Bbb C\}$ is a 
one-dimensional reducing subspace for $L$ and 
we can consider the restriction $L_0$ of $L$ 
to the subspace
$$
H_0=\{x\in L^2(M,\omega): x\perp \bold 1\}=
\{x\in M: \omega(x)=0\}.  
$$
We will show that 
$$
\lim_{t\to\infty}\|\exp{tL_0}\|=0, \tag{4.5}
$$
$\|\cdot\|$ denoting the operator norm in 
$\Cal B(H_0)$.  

Notice that (4.5) implies that $P$ is pure with 
absorbing state $\omega$.  Indeed, for any 
$x\in M$ we set $x_0=x-\omega(x)\bold 1$.  
Then $x_0\in H_0$ and we 
may conclude from (4.5) that 
$$
\lim_{t\to\infty}P_t(x_0)=0,
$$
hence
$$
\lim_{t\to\infty}P_t(x)=\omega(x)\bold 1,
$$
and finally 
$$
\lim_{t\to\infty}\|\rho\circ P_t-\omega\|=0
$$
for every state $\rho$ of $M$ because $M$ is 
finite dimensional.  

In order to prove (4.5), we note first that 
$\{\exp{tL_0}: t\geq 0\}$ is a contraction semigroup 
acting on $H_0$.  Indeed, $\exp{tL}$ is a contraction
in $\Cal B(L^2(M, \omega))$ for every $t$ by virtue
of the inequality 
$$
\|P_t(x)\|^2_{L^2(M,\omega)}=\omega(P_t(x)^*P_t(x)) \leq
\omega(P_t(x^*x))=\omega(x^*x)=\|x\|^2_{L^2(M,\omega)},
$$
and the restriction of $P_t$ to $H_0$ is $\exp{tL_0}$.  

In particular, the spectrum of $L_0$ is contained in 
the left half plane 
$$
\sigma(L_0)\subseteq \{z\in\Bbb C: z+\bar z\leq 0\}.  
$$
We claim that $\sigma(L_0)$ contains no points on 
the imaginary axis $\{iy: y\in \Bbb R\}$.  To see
this, notice first that $0\notin\sigma(L_0)$.  
Indeed, if $L(x)=L_0(x)=0$ for $x\in H_0$ then
$x$ must be a scalar multiple of $\bold 1$ by 
ergodicity, and since $\omega(x)=0$ we have 
$x=0$.  

Suppose now that $\alpha$ is a nonzero real 
number such that $i\alpha\in\sigma(L_0)$.  
Then there is an element $x\neq 0$ in $H_0$ 
for which $L(x)=i\alpha x$.  Note first that 
$x$ is a scalar multiple of a unitary operator.  
Indeed, from the equation $L(x)=i\alpha x$ it 
follows that 
$$
P_t(x)=e^{i\alpha t}x \qquad \text{for every }t\geq 0,
$$
hence 
$$
x^*x=P_t(x)^*P_t(x)\leq P_t(x^*x)
$$
by the Schwarz inequality.  Since 
$\omega(P_t(x^*x)-x^*x)=0$ and $\omega$ is faithful
we conclude that $P_t(x^*x)=x^*x$; so by ergodicity 
$x^*x$ must be a scalar multiple of $\bold 1$.  Thus
$x$ must be proportional to an isometry in $M$.  

We have located a unitary operator $u\in M$ such 
that $L_0(u)=i\alpha u$.  Now we assert that 
$u$ must commute with the self-adjoint set of 
operators $\{v_1,\dots,v_r,v_1^*,\dots,v_r^*\}$.  
To see that we make use of the formula
$$
L(xx^*)-xL(x)^*-L(x)x^*+xL(\bold 1)x^* =
\sum_{j=1}^r[v_j,x][v_j,x]^* \tag{4.6}
$$
(see the proof of Proposition 4.3).  Setting 
$x=u$ we find that the left side of (4.6) is 
$$
-uL(u)^*-L(u)u^* = i\alpha\bold 1-i\alpha\bold 1=0,
$$
and hence 
$$
\sum_{j=1}^r[v_j,u][v_j,u]^* = 0,
$$
from which we deduce that $[v_j,u]=0$ for every $k$.  
Since $u$ is unitary the assertion follows.  

Set 
$$
h=1/2\sum_{j=1}^rv_jv_j^*.  
$$
Since $L(\bold 1)=0$ it follows that $k$ has Cartesian 
decomposition of the form $k=-h+\ell$ where 
$\ell^*=-\ell$, hence $L$ decomposes into a sum of the 
form 
$$
L(x)=L_0(x)+[\ell,x]
$$
where 
$$
L_0(x)=\sum_{j=1}^r v_jxv_j^* -hx-xh.  
$$
By what we have just proved, $L_0(u)=uL_0(\bold 1)=0$.   
It follows that the equation $L(u)=i\alpha u$ reduces to 
$$
[\ell,u]=i\alpha u.  \tag{4.7}
$$
Now since $\ell$ is skew-adjoint, 
$v_s=e^{s\ell}$ defines a one-parameter group 
of unitary operators in $M$ and (4.7) implies that 
for every $s\in\Bbb R$ we have 
$$
v_suv_s^*=e^{i\alpha s}u.  
$$
Since $x\mapsto v_sxv_s^*$ is a $*$-automorphism of 
$M$ for every $s\in\Bbb R$ it follows that the 
spectrum of $u$ must be invariant under all rotations
of the unit circle of the form 
$\lambda\mapsto e^{i\alpha s}\lambda$, contradicting 
the fact that the spectrum of an $N\times N$ unitary 
matrix is a finite subset of 
$\{\lambda\in\Bbb C:|\lambda|=1\}$.  This contradiction
shows that $\sigma(L_0)$ cannot meet the imaginary axis.  

We conclude that 
$$
\sigma(L_0)\subseteq \{z\in\Bbb C: z+\bar z<0\}
$$
and hence there is a positive number $\epsilon$ such 
that 
$$
\sigma(L_0)\subseteq\{z\in \Bbb C: z+\bar z<-2\epsilon\}.  
\tag{4.8}
$$
Consider the operator $A=\exp{L_0}\in\Cal B(H_0)$.  
By the spectral mapping theorem the spectral radius 
of $A$ satisfies
$$
\sup\{|e^z|:z\in\sigma(L_0)\}<e^{-\epsilon}
$$
and hence there is a constant $c>0$ such that 
$$
\|A^n\|\leq c e^{-n\epsilon}, \qquad n=0,1,2,\dots.  
$$
Letting $[t]$ denote the greatest integer not exceeding
$t\geq 0$ we find that for every $t>0$
$$
\|\exp{tL_0}\|\leq \|\exp{[t]L_0}\|=\|A^{[t]}\|\leq 
c e^{-[t]\epsilon},
$$
and hence 
$$
\lim_{t\to\infty}\|\exp{tL_0}\|=0,
$$
as asserted.\qed
\enddemo

\subheading{5.  Applications}

In \cite{4}, a numerical index $d_*(P)$ was introduced
for arbitrary CP semigroups $P=\{P_t: t\geq 0\}$ 
acting on $\Cal B(H)$.  It was 
shown that for unital CP semigroups $P$,
$d_*(P)$ is a nonnegative integer or 
$\infty=\aleph_0$, or $2^{\aleph_0}$, 
and in fact $d_*(P)$ agrees with 
the index of the minimal dilation of $P$ to an \esg.  
In \cite{5}, $d_*(P)$ is calculated in all cases 
where the generator of $P$ is bounded, and in particular 
for CP semigroups acting on matrix algebras.  

We will make use of this numerical index in the following 
result, from which we will deduce Theorem A.  

\proclaim{Theorem 5.1}
Let $\omega$ be a faithful state of $M_r(\Bbb C)$, 
$r\geq 2$, and let $n$ be a positive integer satisfying 
$n\leq r^2-1$.  Then there is a pure CP semigroup 
$P=\{P_t: t\geq 0\}$ acting on $M_r(\Bbb C)$ satisfying 
\roster
\item"{(i)}"
$\omega\circ P_t=\omega$ for every $t\geq 0$, and 
\item"{(ii)}"
$d_*(P)=n$.  
\endroster
\endproclaim

We have based the proof of Theorem 5.1 on the following 
result.  

\proclaim{Proposition 5.2}
Suppose that $T$ is a non-scalar matrix in $M_r(\Bbb C)$, 
$r\geq 2$, and let $\lambda=e^{2\pi i/r}$.  Then there 
is a pair $u,v$ of unitary operators in $M_r(\Bbb C)$ with
the properties
\roster
\item"{5.2.1}"
$u^r = v^r = \bold 1$, 
\item"{5.2.2}"
$vu = \lambda uv$
\item"{5.2.3}"
$\{T,u\}^\prime = \Bbb C\cdot \bold 1$.  
\endroster
\endproclaim

\demo{proof of Proposition 5.2}  The assertion 5.2.3 is that 
the only operators commuting with both $u$ and $T$ are scalars.  
Let $H$ be an $r$-dimensional Hilbert space and identify 
$M_r(\Bbb C)$ with $\Cal B(H)$.  

We claim first that there is an orthonormal basis 
$\xi_0,\xi_1,\dots,\xi_{r-1}$ for $H$ such that 
$$
\<T\xi_0,\xi_k\> \neq 0, \qquad 1\leq k\leq r-1.  \tag{5.3}
$$
Indeed, since $T$ is not a scalar there must be a unit vector 
$\xi_0\in H$ which is not an eigenvector of $T$.  
Thus there is a complex number $a$ and a nonzero 
vector $\zeta$ orthogonal to $\xi_0$ such that 
$$
T\xi_0=a\xi_0 + \zeta.  
$$
Let $c_1,c_2,\dots,c_{r-1}$ be any sequence of nonzero 
complex numbers satisfying 
$$
|c_1|^2+|c_2|^2+\dots+|c_{r-1}|^2=\|\zeta\|^2.  
$$
Since $\zeta\neq 0$ we can find an orthonormal basis 
$\xi_1,\xi_2,\dots,\xi_{r-1}$ for $[\xi_0]^\perp$ such 
that $\<\zeta,\xi_k\>=c_k$ for $k=1,2,\dots,r-1$.  For 
such a choice, 
the set $\{\xi_0,\xi_1,\dots,\xi_{r-1}\}$ is an 
orthonormal basis with the asserted property (5.3).  

Now define $u,v\in\Cal B(H)$ by
$$
\align
u\xi_k&=\lambda^{-k}\xi_k \qquad\text{and}\\
v\xi_k&=\xi_{k\dotplus 1}
\endalign
$$
for $0\leq k\leq r-1$, where $\dotplus$ denotes addition
modulo $r$.  It is obvious that $u$ and $v$ are unitary 
operators, and a straightforward computation shows that 
they satisfy formulas 5.2.1 and 5.2.2.  

We claim now that if $B\in\Cal B(H)$ satisfies 
$BT=TB$ and $Bu=uB$ then $B$ must be a scalar multiple 
of the identity.  Indeed, from $Bu=uB$ and the fact 
that $u$ is a unitary operator with distinct 
eigenvalues, we find that 
each $\xi_k$ must be an eigenvector of 
both $B$ and $B^*$.  Choosing
$d_k\in\Bbb C$ such that $B\xi_k=d_k\xi_k$, 
then $B^*\xi_k=\bar d_k\xi_k$ and 
for each $k=1,2,\dots,r-1$ we have 
$$
d_0\<T\xi_0,\xi_k\>=\<TB\xi_0,\xi_k\>=
\<BT\xi_0,\xi_k\>=\<T\xi_0,B^*\xi_k\>=d_k\<T\xi_0,\xi_k\>.  
$$
It follows that $(d_k-d_0)\<T\xi_0,\xi_k\>=0$ 
for $1\leq k\leq r-1$.  Because none of the inner products
$\<T\xi_0,\xi_k\>$ can 
be zero we conclude that $d_0=d_1=\dots=d_{r-1}$.  
Thus $B=d_0\cdot\bold 1$, establishing proposition 
5.2.  \qed
\enddemo

\remark{Remarks}
Let $\lambda$ be a primitive $r$th root of unity and 
let $u,v$ be two unitaries satisfying condition 5.2.1 
and 5.2.2.  Consider the family of $r^2$ 
unitary operators $\{w_{i,j}: 0\leq i,j\leq r-1\}$
defined by
$$
w_{i,j}=u^iv^j.  
$$
We may consider that the indices $i,j$ range over the 
abelian group $\Bbb Z/r\Bbb Z$, and with that convention
the $w_{i,j}$ are seen to 
satisfy the commutation relations for this 
group
$$
\align
w_{i,j}w_{p,q}&=\lambda^{jp}w_{i+ p,j+q}\tag{5.4}\\
w_{i,j}^* &= \lambda^{ij}w_{-i,-j}\tag{5.5}
\endalign
$$
where the operations $i+p$, $j+q$, $-i$, $-j$ 
are performed modulo $r$.  
Of course, we have $w_{0,0}=\bold 1$.  It follows from (5.4) 
and (5.5) that the set of operators $\{w_{i,j}\}$ satisfies 
$$
w_{i,j}w_{p,q}w_{i,j}^* = \lambda^{jp-qi}w_{p,q}.  
$$
This formula, together with the fact that $\lambda$ is 
a primitive $r$th root of unity, implies that 
$$
\text{trace}(w_{p,q}) = 0, \qquad \text{for }0\leq p,q\leq r-1,
\quad p+q>0.  \tag{5.6}
$$
In particular, from (5.4)--(5.6) we see that relative to the 
inner product on $M_r(\Bbb C)$ defined by the normalized
trace, the set of operators $\{w_{i,j}: 0\leq i,j\leq r-1\}$ is 
an orthonormal basis.  Thus the $\{w_{i,j}: 0\leq i,j\leq r-1\}$
are linearly independent.  
\endremark

\demo{proof of Theorem 5.1}
Assume first that $\omega$ is not the tracial state, and 
let $\Omega$ be its density matrix.  Then $\Omega$ is not 
a scalar multiple of the identity and Proposition 5.2 provides 
a pair of unitary operators $u,v$ satisfying (5.2.1), (5.2.2) 
and (5.2.3) for $T=\Omega$.  
Define $w_{i,j}=u^iv^j$, $0\leq i,j\leq r-1$.  
By the preceding remarks the set of $r^2-1$ unitary operators 
$\Cal S=\{w_{i,j}: 0\leq i,j\leq r-1, i+j>0\}$ 
is linearly independent
and consists of trace zero operators.  

Choose $n$ satisfying $1\leq n\leq r^2-1$ and let 
$v_1, v_2, \dots,v_n$ be any set of $n$ distinct 
elements of $\Cal S$ such that $v_1=w_{1,0}=u$.  By 
(5.2.3) we have 
$$
\{\Omega,v_1\}^\prime = \Bbb C \bold 1,   
$$
and hence 
$$
\{\Omega,v_1,v_2,\dots,v_n\}^\prime = \Bbb C \bold 1. \tag{5.7}
$$

Consider the completely positive map of $M_r(\Bbb C)$ defined 
by
$$
Q(x) = \Omega^{-1/2}(\sum_{k=1}^n v_k x v_k^*)\Omega^{-1/2}.  
$$
Since the $v_k$ are unitary operators we have
$$
\sum_{k=1}^r v_kv_k^* = \sum_{k=1}^r v_k^*v_k,
$$
hence Corollary 3.16 implies
that there is an operator $k\in M_r(\Bbb C)$ such that 
$$
L(x) = Q(x) + kx + xk^*
$$
generates a unital CP semigroup $P=\{P_t: t\geq 0\}$ satisfying 
$\omega\circ P_t = \omega$ for every $t\geq 0$.  Because of 
(5.7), Theorem 4.4 implies that $P$ is a pure semigroup.  

It remains to show that $d_*(P) =n$, and for that we appeal to the
results of \cite{5}.  Consider the linear span
$$
\Cal E = 
\text{span}\{\Omega^{-1/2}v_1,\Omega^{-1/2}v_2,\dots,\Omega^{-1/2}v_n\}.
$$
We claim first that 
$\Cal E\cap\Bbb C\bold 1 = \{0\}$.  Indeed, if this intersection 
were not trivial then we would have 
$$
\bold 1 = c_1\Omega^{-1/2}v_1+\dots+c_n\Omega^{-1/2}v_n
$$
for some scalars $c_1,\dots,c_n$.  Hence 
$$
\Omega^{1/2}=c_1v_1+\dots+c_nv_n.  
$$
This is impossible because the left side has positive trace, 
while by (5.6) the right side has trace zero.

We can make $\Cal E$ 
into a metric
operator space \cite{4, Definition 1.9} 
by declaring the linear basis 
$\Omega^{-1/2}v_1,\dots,\Omega^{-1/2}v_n$ 
to be an orthonormal basis, 
and once this is done we find that $\Cal E$ is the metric operator
space associated with the completely positive map $Q$.  From 
\cite{5, Theorem 2.3} we have $d_*(P) = \dim \Cal E = n$, 
as required.  

It remains to deal with the case where $\omega$ is the normalized 
trace on $M_r(\Bbb C)$.  That requires a small variation 
of the preceding argument.  Choose an arbitrary operator 
$T\in M_r(\Bbb C)$ so that $T$ is not a scalar and satisfies 
$T^*=-T$.  Let $\lambda$ be a primitive $r$th root of 
unity and let $u$, $v$ be two unitary operators satisfying 
the three conditions of Proposition 5.2.  Now we form the 
operators $w_{i,j}$ exactly as before, and obtain $n$ 
unitary operators $\{v_1,v_2,\dots,v_n\}$ by enumerating 
the elements of $\{w_{i,j}: 0\leq i,j\leq r-1, i+j>0\}$
in such a way that $v_1=u$.  Define an operator $L$ on 
$M_r(\Bbb C)$ by
$$
L(x)=\sum_{k=1}^n v_kxv_k^* -nx + [T,x].  
$$
Notice that $L(\bold 1) =0$ and, since we obviously have 
$\sum_k v_kv_k^* = \sum_k v_k^*v_k$, it follows that 
$\text{trace}(L(x))=0$ for all $x\in M_r(\Bbb C)$.  Hence 
$L$ is the generator of a unital 
CP semigroup $P=\{P_t: t\geq 0\}$ which preserves the 
tracial state $\omega$.  

Notice that $P$ is pure.  Indeed, by (5.2.3) we have 
$\{v_1,T\}^\prime = \Bbb C\bold 1$, and hence the 
$*$-algebra generated by the set $\{v_1,\dots,v_n,T\}$ is 
irreducible.  Theorem 4.4 implies that $P$ is a pure CP semigroup.  

Finally, $d_*(P)=n$ follows exactly as in the non-tracial case
already established. \qed
\enddemo

We are now in position to prove Theorem A, as stated in 
the introduction.  Let $r$ and $n$ be 
positive numbers with $r\geq 2$, and let 
$\lambda_1,\lambda_2,\dots,\lambda_r$ be a sequence of positive 
numbers summing to $1$.  We have to show that there is a cocycle 
perturbation of the $CAR/CCR$ flow of index $n$ which has an 
absorbing state with eigenvalue list 
$\lambda_1,\lambda_2,\dots,\lambda_r$.  

We first consider the case in which $n\leq r^2-1$.  Let 
$H_0$ be a Hilbert space of dimension $r$, and identify 
$M_r(\Bbb C)$ with $\Cal B(H_0)$.  Choose an orthonormal basis 
$\xi_1,\xi_2,\dots,\xi_r$ for $H_0$ and let $\omega_0$ be the 
state of $\Cal B(H_0)$ defined by 
$$
\omega_0(x)=\sum_{k=1}^r \lambda_k\<x\xi_k,\xi_k\>.  
$$
Then $\omega_0$ is a faithful state 
on $\Cal B(H_0)$ having eigenvalue list
$\lambda_1,\lambda_2,\dots,\lambda_r$.  By Theorem 
5.1, there is a pure CP semigroup $P=\{P_t: t\geq 0\}$ acting 
on $\Cal B(H_0)$ such that $\omega_0\circ P_t=\omega_0$ for 
every $t\geq 0$.  Using Bhat's dilation theorem \cite{7,8}, 
there is a Hilbert space $H\supseteq H_0$ and an \esg\ 
$\alpha=\{\alpha_t: t\geq 0\}$ acting 
on $\Cal B(H)$ such that if we identify $\Cal B(H_0)$ with 
the corner $p_0\Cal B(H)p_0$ ($p_0$ denoting the projection 
of $H$ onto $H_0$), then we have $\alpha_t(p_0)\geq p_0$ for 
every $t\geq 0$ and for every $x\in\Cal B(H_0)$ 
$$
P_t(x)=p_0\alpha_t(x)p_0, \qquad t\geq 0.  
$$
Using \cite{2}, we may assume that $\alpha$ is 
{\it minimal} over the projection $p_0$.  

Now by Proposition 2.4, $\alpha$ is a pure \esg.  Moreover, 
if we define a normal state $\omega$ of $\Cal B(H)$ by 
$$
\omega(x) = \omega_0(p_0xp_0), 
$$
then $\omega$ must be invariant 
under $\alpha$.  Indeed, since 
$\alpha_t(p_0)\geq p_0$ we have for every $x\in\Cal B(H)$
$$
p_0\alpha_t(x)p_0=p_0\alpha_t(p_0xp_0)p_0 = P_t(p_0xp_0),
$$
hence 
$$
\omega(\alpha_t(x))=\omega_0(P_t(p_0xp_0)=\omega_0(p_0xp_0)
=\omega(x), 
$$
as asserted.  By the general discussion of 
section 1 it follows that $\omega$ is an absorbing
state, and of course the eigenvalue list of $\omega$ is the 
same as that for $\omega_0$, namely 
$\lambda_1,\lambda_2,\dots,\lambda_r$.  Thus it only 
remains to show that $\alpha$ is conjugate to a cocycle 
perturbation of the $CAR/CCR$ flow of index $n$.  But by 
Corollary 4.21 of \cite{5}, $\alpha$ is cocycle
conjugate to a $CAR/CCR$ flow of index $d_*(P)=n$, and 
the proof of this case is complete.  

Suppose now that $n>r^2-1$.  In this case, pick any positive
integer $k\leq r^2-1$.  By what was just proved,
we can find a cocycle perturbation $\alpha$ of the $CAR/CCR$ 
flow of index $k$ which has an absorbing state $\omega$ having 
eigenvalue list $\lambda_1,\lambda_2,\dots,\lambda_r$.  Moreover,
letting $p_0$ be the support projection of $\omega$ then $p_0$ 
has rank $r$ and if $P$ is the CP semigroup obtaind by 
compressing $\alpha$ to $p_0\Cal B(H)p_0$, then $P$ is 
a pure CP semigroup and $\alpha$ can be assumed to be 
the minimal dilation of $P$.  

We will show how to use $\alpha$ construct a 
{\it nonminimal} dilation $\beta$ of $P$ which is pure, 
conjugate to a cocycle perturbation of the $CAR/CCR$ flow
of index $n$, and has an absorbing state with the same 
eigenvalue list.  For that, let $m=n-k$ and let $\alpha^m$
be the $CAR/CCR$ flow of index $m$, acting on $\Cal B(K)$.  
It is known that 
every $CAR/CCR$ flow has a pure absorbing state $\rho$
(the vacuum state) \cite{13}.  
Thus letting $\zeta\in K$ be 
the vacuum vector then we have 
$$
\rho(x)=\<x\zeta,\zeta\>.  
$$ 
If we write $[\zeta]$ for the rank-one projection defined
by $\zeta$ then $\alpha^m_t([\zeta])\geq[\zeta]$
for every $t\geq 0$ and in fact 
$$
\lim_{t\to\infty}\alpha^m_t([\zeta])=\bold 1_K. \tag{5.8}
$$  

Let $\beta$ be the \esg\ defined on $\Cal B(H\otimes K)$ 
by $\beta=\alpha\otimes \alpha^m$, i.e., 
$$
\beta_t(x\otimes y)=\alpha_t(x)\otimes\alpha^m_t(y), 
\qquad x\in\Cal B(H), y\in\Cal B(K), t\geq 0.  
$$
$\beta$ is obviously a 
cocycle perturbation of the $CAR/CCR$ flow of
index $n=k+m$.  We will show that $\beta$ is a pure
\esg\ having an invariant state with eigenvalue list
$\lambda_1,\lambda_2,\dots,\lambda_r$.  

To that end, consider the normal state $\omega^\prime$ 
defined on $\Cal B(H\otimes K)$ by 
$$
\omega^\prime = \omega\otimes \rho.  
$$
Since $\rho$ is a vector state, $\omega^\prime$ has the 
same eigenvalue list as $\omega$, namely 
$\lambda_1,\lambda_2,\dots,\lambda_r$.  Moreover, 
$\omega^\prime$ is invariant under $\beta$ because 
$\omega$ (resp. $\rho$) is invariant under 
$\alpha$ (resp. $\alpha^m$).  
Thus it remains to show that $\beta$ is a pure 
\esg .  

For that, we appeal to Proposition 2.4 as follows.  
Let $q_0=p_0\otimes[v]$ be the support projection 
of $\omega^\prime$.  Then we have 
$$
\beta_t(q_0)=\alpha_t(p_0)\otimes\alpha^m_t([v]).  
$$
Since the projections $\alpha_t(p_0)$ (resp. 
$\alpha^m_t([v])$) increase with $t$ to $\bold 1_H$ 
(resp. $\bold 1_K$), it follows that 
$\beta_t(q_0)\geq q_0$ and 
$$
\lim_{t\to\infty}\beta_t(q_0)=\bold 1_{H\otimes K}.  
$$
Thus if we let $Q=\{Q_t: t\geq 0\}$ be the CP semigroup 
obtained by compressing $\beta$ to the corner 
$q_0\Cal B(H\otimes K)q_0$, it follows that $\beta$ 
is a (nonminimal) dilation of $Q$.  Finally, since 
$[v]$ is one-dimensional, $Q$ is conjugate to the original 
CP semigroup $P$, and is therefore pure.  By 
Proposition 2.4, we conclude that $\beta$ is a 
pure \esg.

We have established all but the third paragraph of 
Theorem A, to which we now turn our attention.  Let 
$r\geq 2$ be an integer and let $\beta$ be an \esg\ 
acting on $\Cal B(H)$, $H$ being a separable infinite 
dimensional Hilbert space, which has 
an absorbing state with eigenvalue list 
$\lambda_1,\lambda_2,\dots,\lambda_r$.  Assuming 
that $\beta$ is minimal over the support projection 
$p_0$ of $\omega$, we have to show that $\beta$ is 
cocycle conjugate to a $CAR/CCR$ flow of index $n$ 
where $n$ is a positive integer not exceeding 
$r^2-1$.  

Let $H_0=p_0H$ and let 
$P=\{P_t: t\geq 0\}$ be the CP semigroup obtained
by compressing $\beta$ to the corner 
$p_0\Cal B(H)p_0\cong \Cal B(H_0)$.  
Let $L$ be the generator of the semigroup $P$.  By 
\cite{4,5} there is an operator $k\in \Cal B(H_0)$ 
and a metric operator space 
$\Cal E\subseteq \Cal B(H_0)$ (possibly $\{0\}$) satisfying 
$\Cal E\cap \Bbb C\bold 1=\{0\}$ and which give rise to
$L$ as follows:
$$
L(x) = \sum_{k=1}^nv_kxv_k^* + kx + xk^*, 
\qquad x\in \Cal B(H_0),  \tag{5.9}
$$
$v_1,v_2,\dots,v_n$ denoting any orthonormal basis 
for $\Cal E$.  Since $\Cal E$ is 
a proper subspace of the $r^2$-dimensional vector
space $\Cal B(H_0)$, the integer $n=\dim\Cal E$ 
has possible values $0,1,\dots,r^2-1$.  

Note first that $n$ cannot be $0$.  For in that 
case (5.9) reduces to $L(x)=kx+xk^*$.  Using the fact
that $L(\bold 1) =0$, we find that $k$ must be a skew-adjoint 
operator for which $L(x)=[k,x]$, hence
$$
P_t(x)=\exp{tL}(x)=e^{tk}xe^{-tk} 
$$
is a semigroup of $*$-automorphisms of 
$\Cal B(H_0)$.  Since $\beta$ is a minimal 
dilation of $P$ we must have $H=H_0$ and 
$\beta_t=P_t$ for every $t\geq 0$, contradicting the 
fact that $\beta$ is an \esg\ acting 
on an infinite dimensional type $I$ factor.  

Thus $1\leq n\leq r^2-1$.  Theorem 2.3 of \cite{5}
implies that the index of $P$ is given by 
$d_*(P)=\dim\Cal E=n$, and by \cite{4} Theorem 4.9 we 
have $d_*(\beta) = d_*(P)=n$.   $\beta$ must be completely 
spatial by \cite{5} Theorem 4.8, and finally 
by the classification results of \cite{1} 
(Corollary of Proposition 7.2) every 
completely spatial \esg\ is conjugate to a cocycle 
perturbation of a $CAR/CCR$ flow.  That completes the 
proof of Theorem A.

%\vfill
%\pagebreak

\end